\documentclass[10pt, english]{amsart}

\usepackage{amsmath,amssymb,enumerate}

\usepackage[T1]{fontenc}
\usepackage[all,cmtip]{xy}

\usepackage{babel}
\usepackage{amstext}
\usepackage{amsmath}
\usepackage{amsfonts}
\usepackage{latexsym}
\usepackage{ifthen}

\usepackage{xypic}
\xyoption{all}
\newcommand{\sJ}{{\mathcal J}}
\newcommand{\sK}{{\mathcal K}}

\newtheorem{lemma1}{}[section]

\newenvironment{theorem}{\begin{lemma1}{\bf Theorem.}}{\end{lemma1}}

\newenvironment{remark}{\begin{lemma1}{\bf Remark.}\rm}{\end{lemma1}}

\newenvironment{the local obstruction - setup}{\begin{lemma1}{\bf The local obstruction - setup.}}{\end{lemma1}}

\newenvironment{remark*}{{\bf Remark.}}{}
\newenvironment{remarks*}{{\bf Remarks.}}{}
\newenvironment{example*}{{\bf Example.}}{}
\newenvironment{assumption*}{{\bf Assumption.}}{}

\newcommand{\Q}{\ensuremath{\mathbb{Q}}}

\newcommand{\C}{\ensuremath{\mathbb{C}}}
\newcommand{\N}{\ensuremath{\mathbb{N}}}
\newcommand{\PP}{\ensuremath{\mathbb{P}}}

\usepackage{eurosym}

\newcommand{\holom}[3]{\ensuremath{#1:#2  \rightarrow #3}}
\newcommand{\fibre}[2]{\ensuremath{#1^{-1} (#2)}}

\makeatletter
\ifnum\@ptsize=0  \fi \ifnum\@ptsize=2  \fi 

\newcommand\sF{{\mathcal F}}
\newcommand\sG{{\mathcal G}}

\newcommand\sI{{\mathcal I}}

\newcommand\sO{{\mathcal O}}

\usepackage{xcolor}
\usepackage{hyperref}

\author{Andreas H\"oring}
\author{Thomas Peternell}

\address{Andreas H\"oring, Universit\'e C\^ote d'Azur, CNRS, LJAD, France}
\email{Andreas.Hoering@univ-cotedazur.fr}

\address{Thomas Peternell, Mathematisches Institut, Universit\"at Bayreuth, 95440 Bayreuth, 
Germany}
\email{thomas.peternell@uni-bayreuth.de}

\title{A contraction theorem for divisors fibering over a curve} 
\date{September 24, 2024}

\begin{document}

\maketitle

\begin{abstract}
Given a $\Q$-Cartier divisor $S \subset X$ admitting a fibration $S \rightarrow B$
onto a curve we give sufficient conditions for the existence of a bimeromorphic
contraction $X \rightarrow Y$ contracting $S$ onto $B$. As a corollary we recover
a contraction result for compact K\"ahler threefolds.
\end{abstract}

\section{Introduction} 

In this paper we give a proof of the following theorem which is part of the minimal model program for K\"ahler threefolds.

\begin{theorem} \label{thmA}  Let $X$ be a $\mathbb Q$-factorial normal compact K\"ahler threefold with at most terminal 
singularities. Let $R$ be a divisorial $K_X$-negative extremal ray on $X$, and assume that the divisor $S \subset X$ covered by curves $C \subset X$ such that  $[C] \in R$ admits a fibration 
$$
f: S\to B
$$ 
onto a curve contracting exactly these curves.

Then there exists a normal compact threefold $Y$ and a bimeromorphic map $\varphi: X \to Y$ such that $\varphi_{\vert S}$ is the Stein factorisation of $f$ and 
$\varphi_{\vert X \setminus S}$ 
is biholomorphic. 
\end{theorem}

This result was initially claimed in \cite[Cor.7.9]{HP16} and deduced from \cite[Prop. 7.4]{HP16} whose proof is based on \cite[Thm.2]{AV84}.
It has been pointed out to us by Matei Toma that \cite[Thm.2]{AV84} and \cite[Prop.7.4]{HP16} are false as stated, as shown by a surprising example
of Fujiki \cite[Prop.3]{Fuj75}. One observes that Fujiki's example needs a cohomological assumption that is typically not satisfied for contractions appearing in the minimal model program (cf. Remark \ref{remark-fujiki}). Using quite sophisticated techniques from MMP, Das and Hacon \cite{DH23} give a proof of the contraction theorem for K\"ahler threefolds, thereby repairing this gap.

In this paper we follow a different strategy which is not limited to the threefold setting and might be of general interest: the bimeromorphic contraction $\varphi$ always exists topologically. Assuming that $\varphi$ is holomorphic in the complement of finitely many fibres of $S \rightarrow B$ (a condition that is easily checked for terminal threefolds), we 
aim to extend the complex structure to get a holomorphic map. We achieve this goal by extending the morphism $f$
to a map between infinitesimal neighbourhoods $f_k: S_k \rightarrow B_k$ for arbitrary $k \in \N$ and conclude with a result
of Bingener \cite{Bin81}.
This strategy applies whenever the ``exceptional'' divisor fibres over a curve, so we obtain a contraction theorem in arbitrary dimension (see Theorem \ref{thmB}).

{\bf Acknowledgements.} We thank Matei Toma for making us aware of Fujiki's classical example.

\section{Main result and its proof}

\begin{theorem} \label{thmB}  Let $X$ be a normal compact complex space. Let $S \subset X$ be a prime divisor that admits a morphism 
$$
f: S\to B
$$ 
onto a curve $B$ such that $f_* \sO_S \simeq \sO_B$.  Assume that
\begin{itemize}
\item the divisor $S$ is $\Q$-Cartier of index $m$ and the line bundle $\sO_S(-mS)$ is $f$-ample; and
\item there exists a finite set $Z \subset B$ and 
a bimeromorphic map 
$$
\varphi^\circ: X^\circ := (X \setminus \fibre{f}{Z}) \to Y^\circ
$$
such that $\varphi^\circ=f$ on $S \setminus \fibre{f}{Z}$ and $\varphi^\circ_{\vert X \setminus S}$ 
is biholomorphic; and
\item for every $k \in \N$ the natural map $\sI_{B \setminus Z}^k \rightarrow \varphi^\circ_* (\sI_{S \setminus \fibre{f}{Z}}^k)$ is an isomorphism
in the generic point of $B$.
\end{itemize}
Then there exists a normal compact complex space $Y$ and a bimeromorphic map $\varphi: X \to Y$
such that $\varphi_{\vert S} = f$ and $\varphi_{\vert X \setminus S}$ 
is biholomorphic. 
\end{theorem}

\begin{remark*}
The third condition in Theorem \ref{thmB} holds e.g. if $X^\circ$ is the blow-up of $Y^\circ$ along the ideal sheaf $\sI_B$ \cite[II, \S 3]{AT82}.
\end{remark*}

\begin{proof}
Set
$$
X^{\circ} := X \setminus \fibre{f}{Z}, \ S^{\circ} := S \setminus \fibre{f}{Z}, \ B^\circ := B \setminus Z
$$
and $f^{\circ} = f \vert S^{\circ}$.

We introduce an equivalence relation on $X$ by identifying $x_1$ and $x_2$ if and only if $f(x_1) = f(x_2)$. Let $Y$ be the quotient space with 
quotient map $\varphi$. 
Then $Y$ is a compact Hausdorff space
 In fact, given $y_1, y_2 \in Y$, then there exist open neighborhoods $U_i$ of $\varphi^{-1}(y_i)$ which are disjoint. 
Further,  $Y$ contains $Y^{\circ} $ as open set, and $\varphi_{\vert X^{\circ}} = \varphi^{\circ}$. 
Setting $\sO_Y = \varphi_*(\sO_X)$, we obtain a locally ringed space $(Y,\sO_Y)$. Note that it is a priori not at all clear whether $\sO_Y$ is coherent as $\sO_Y$-module.

Let $\sI_S$ be the  ideal sheaf of $S \subset X$, giving rise to infinitesimal neighborhoods
$$ S_k := (S,\sO_X/\sI_S^k)).$$
We further set
$$ \sI_B := \varphi_*(\sI_S)$$ 
and introduce the ringed spaces
$$ B_k := (B,\sO_Y/\sI_B^k).$$

{\bf Step 1.} 
Having in mind that $B$ is a possibly singular compact complex curve, we show 
$$ \sO_B = \sO_{B_1},$$ 
hence $B_1$ is a complex space.
Starting with the exact sequence
$$ 0 \to \sI_S \to \sO_X \to \sO_S \to 0$$
and applying $\varphi_*$, we obtain the exact sequence
$$ 0 \to \sI_B = \varphi_*(\sI_S) \to \sO_Y = \varphi_*(\sO_X) \buildrel {\nu} \over {\to } \sO_B = \varphi_*(\sO_S )\to R^1 \varphi_*(\sI_S).$$ 
Hence we need to show that $\nu$ is surjective. But this is clear since the $\sO_B$-module ${\rm Im} \  \nu$ contains the constant section $1_B$, hence ${\rm Im}\  \nu = \sO_B$.

We have a canonical injective map for all positive integers $k$, 
$$ \sI_B^k \to \varphi_*(\sI_S^k)$$ 
with finitely supported cokernel (supported on $Z$). 

This gives rise to a canonical map
$$ \mu_k: \sI_B^k/\sI_B^{k+1} \to \varphi_*(\sI_S^k/\sI_S^{k+1})$$
which is an isomorphism over $B^{\circ}$. 
The
existence of $\mu_k$  follows from the canonical morphism  $\varphi_* \sI_S^k/\varphi_* \sI_S^{k+1} \rightarrow \varphi_* (\sI_S^k/\sI_S^{k+1})$.

Set 
$$ Q:= {\rm Coker} \mu_1, \qquad R := \ker \mu_1.$$

{\bf Step 2.} 
We claim that $\sI_B^k/\sI_B^{k+1}$ is a coherent $\sO_B$-module for all $k \in \N$. 

{\bf Step 2a.} 
We first reduce to the case  $ k =1$. In fact, there is a canonical epimorphism
$$  \lambda_k: S^k(\sI_B/\sI_B^2) \to \sI_B^k/\sI_B^{k+1}, $$
whose kernel $\mathcal K_k$ is supported on $Z$. Since $\sI_B/\sI_B^2$ is coherent, so is $S^k(\sI_B/\sI_B^2)$. Since $Z$ is finite, $\mathcal K_k$ is clearly 
coherent, hence $\sI_B^k/\sI_B^{k+1}$ is coherent. 
Indeed, on the noetherian scheme $B$ we can use \cite[II, Prop.5.7]{Har77}.

{\bf Step 2b.} 
To show that $\sI_B/\sI_B^2$ is coherent, note that 
since $\varphi_*(\sI_S/\sI_S^{2}) = f_* (\sI_S/\sI_S^{2})$ is a coherent sheaf on $B$, by the usual yoga of coherent sheaves, the claim is equivalent to   $$ h^0(Y, R) < \infty$$ and 
 $$ h^0(Y, Q) < \infty.$$
This is to say that at every point $b \in Z$, the stalk $R_b$ resp. $Q_b$  is a finite-dimensional $\mathbb C$-vector space. 
Since $Q$ is a quotient of $\varphi_*(\sI_S/\sI_S^2)  $, the assertion for $Q$ is clear and it remains to treat $R$. 

First notice that 
$$ R \simeq \varphi_*(\sI_S^2)/\varphi_*(\sI_S)^2.$$ 
Consider the canonical morphisms
$$ \tau_1:  \varphi^* \varphi_*(\sI_S) \to \sI_S$$ 
and
$$ \tau_2:  \varphi^* \varphi_*(\sI_S^2) \to \sI_S^2$$ 
which are surjective on $X^{\circ}$. 

{\it Claim.}  There is an open neighborhood $U $ of $b \in B$ and a coherent sheaf $\sF$ on $\tilde U = \varphi^{-1}(U)$ such that 
$\sF \subset {\rm Im} \tau_2$ and $\sF = \sI_S^2$ on $\tilde U \cap X^{\circ}$. 

Given the claim, it suffices to show
that $\varphi_*(\sI_S^2) /\varphi_*(\sF)$ is coherent, i.e.,
$$h^0(\varphi_*(\sI_S^2) /\varphi_*(\sF)) < \infty.$$ 
 In fact,  $(\varphi_* \sI_S)^2 \subset \varphi_* \sF$ since locally around $b$, the sheaf $(\varphi_*(\sI_S))^2$ is generated by $\varphi_*(h_1) \cdot \varphi_*(h_2)$ with $h_j$ local 
functions around $\varphi^{-1}(b)$. Now use that $\sF \subset {\rm Im}(\tau_2)$. Since $\varphi_*(\sF)$ is coherent and  $\varphi_*(\sI_S^2) /\varphi_*(\sF)$ is supported on $b$, 
 $\varphi_*(\sI_S^2) /\varphi_*(\sF)$is coherent.

Since 
$$ \varphi_*(\sI_S^2) /\varphi_*(\sF) \subset \varphi_*(\sI^2_S/\sF),$$ 
it suffices to show that
$$ h^0(\tilde U,\sI^2_S/\sF) < \infty.$$
But $ \sI^2_S/\sF$ lives on a compact complex space which is set-theoretically contained in $f^{-1}(b)$, so finite-dimensionality follows. 
Alternatively, $\varphi_*(\sI_S^2/\sF) $ is a coherent sheaf on $B$, supported on $Z$. 

It remains to prove the {\it Claim.} 
We have an exact sequence 
$$ 0 \to \ker \mu_1 \to \sI_B/\sI_B^2 \to {\rm Im}  \mu_1 \to 0$$
where ${\rm Im}\  \mu_1 \subset \varphi_*(\sI_S/\sI_S^2)$ is a subsheaf of the coherent sheaf $\varphi_*(\sI_S/\sI_S^2)$  with equality on $B^{\circ}$. 
Since $B$ is noetherian and $\varphi_*(\sI_S/\sI_S^2)$ coherent, the subsheaf ${\rm Im}\  \mu_1$ is also coherent.
Now we take locally near $b$  two general sections of ${\rm Im} \mu_1$, lift them first to sections of $\sI_B/\sI_B^2$ and then to sections 
of $\sI_B$ over some small neighborhood $U$ of $b$. 
These yield sections $s_j$ on $\sI_S$ over $\tilde U$ and 
thus a morphism 
$$ \sO_{\tilde U}^{\oplus 2} \to \sI_S \vert_{\tilde U}$$
whose image is contained in ${\rm Im} \tau_1$ with equality over $\tilde U \cap X^{\circ}$.
Now we let $\sF$ to be the image of the morphism
$$ \sO_{\tilde U}^{\oplus 3} \to \sI_S^2 \vert_{\tilde U}$$
given by $s_1, s_2, s_1 \cdot s_2$;
then clearly the assertions in the Claim are verified and Step 2b is accomplished.

{\bf Step 3.} 
We are going to show that  the ringed space $B_k$ 
is an mFB-space in the sense of Forster-Knorr \cite{FK72} for all $k \in \N$. This notion will be explained in the sequel. 

{\bf Step 3a.} 
First, for any $U \subset B$ open, we have to construct a Fr\'echet topology on $\sO_{B_k}(U)$ for any open set $U \subset B$. 
To  do this, we use the canonical epimorphism
$$  \sO_{B_k} = \sO_Y/\varphi_*(\sI_S)^k \to \sO_Y/\varphi_*(\sI_S^k)$$ 
whose kernel $\sG = \varphi_*(\sI_S^k)/\varphi_*(\sI_S)^k$  is supported on $Z$. Hence $H^1(U,\sG) = 0$ and therefore we obtain an epimorphism
$$u:  \sO_{B_k}(U) \to (\sO_Y/\varphi_*(\sI_S^k))(U).$$ 
By the exact sequence 
$$ 0 \to \varphi_*(\sI_S^k) \to \sO_Y \buildrel {\lambda} \over {\to} \varphi_*(\sO_X/\sI_S^k)  \to R^1\varphi_*(\sI_S^k), $$ 
we see that 
$$ {\rm Im}\lambda = \sO_Y/\varphi_*(\sI_S^k),$$ hence via the above exact sequence
$$  (\sO_Y/\varphi_*(\sI_S^k))(U) \simeq ({\rm Im} \lambda)(U).$$ 
Since  $(\sO_X/\sI_S^k)(f^{-1}(U))$ has a canonical Fr\'echet topology,  see \cite[Thm.5,p.167]{GR79}, so has $ ({\rm Im} \lambda)(U)$ as closed subspace 
{by \cite[p.169]{GR79}}
It suffices to show that $\ker u $ is finite-dimensional,
then $\sO_{B_k}(U)$ will be Fr\'echet, too. 
Now $$\ker F_U = \bigcup_{b \in U \cap Z} \varphi_*{(\sI_S^k)}_b/ {(\sI_B^k)}_b, $$
which is finite-dimensional by Step 2.
Indeed we have an injection 
$$
\varphi_* \sI_S^k/\sI_B^k \hookrightarrow \sI_B/\sI_B^k
$$
and the latter space is a coherent $\sO_B$-module by Step 2 and an induction argument.

Alternatively, by induction, $\sO_{B_k}$ is a coherent $\sO_B$-module by Step1 and 2. Hence $\sO_{B_k}(U)$ carries a unique Fr\'echet structure by \cite[Thm.5,p.167]{GR79}. 

{\bf Step 3b.} 
By definition of an mFB-space,  two conditions have to be verified. First of all, ${(\sO_{B_k})_b}$ is a local ring for all $b \in B$. This is clear. Second, given any $b \in B$, and any open set $U$  containing
$b$ in $B$ there is an open neighborhood $V \subset U$ of $b$ in $B$ such that the restriction map
$$ r: \sO_{B_k}(U) \to \sO_{B_k}(V)$$
has the following property (*). 

{\it For any bounded set $A \subset \sO_{B_k}(U)$ there is an m-bounded set $C \subset \sO_{B_k}(V)$ and a positive number $k$ such that 
$$ r(A) \subset kC.$$}
Recall that a set $C$ is m-bounded, if is bounded and closed, absolutely convex, with $C \cdot C \subset C$ and $C \cap \mathbb C \ne \emptyset$. 

Since complex spaces are mFB-spaces, \cite[p.120]{FK72}, this condition has anyway only to be checked for $b \in Z$. 
We will use the following fact, \cite[p.119]{FK72}. 

{\it Fact.} Suppose $(X,\sO_X)$ is an mFB-space. Then property (*) holds for  {\it any} open sets $V \subset\subset U$, i.e., for any 
bounded set $A \subset \sO_{X}(U)$ there is an m-bounded set $C \subset \sO_{X}(V)$ and a positive number $k$ such that 
$$ r(A) \subset kC.$$

We consider the composed morphism 
$$\tilde u:  \sO_{B_k}(U) \buildrel {u} \over { \to} (\sO_Y/\varphi_*(\sI_S^k))(U) \to    \sO_{S_k}(f^{-1}(U)) $$
Since ${\rm Im} \tilde u $ is closed in $\sO_{S_k}(f^{-1}(U))$ (this is clear, but follows also by a general fact, \cite[p.169]{GR79}, since $\sO_Y/\varphi_*(\sI_S^k)$ is a coherent $\sO_B$-module), 
and since the complex space $S_k$ is an mFB-space,
it has property (*), too, by the above {\it Fact}.
Since $$\sO_{B_k}(U) \simeq \ker \tilde u \oplus {\rm Im} \tilde u $$ as Fr\'echet spaces, and since $\ker \tilde u $ is finite-dimensional, it is easily checked that 
$\sO_{B_k}(U)$ is an mFB-space as well.

{\bf Step 4.} We next show that 
$B_k := (B,\sO_Y/\sI^k)$ is a complex space and that the map $f_k$ of ringed spaces $S_k \to B_k$ is holomorphic.
We argue by induction and will use \cite[10.3]{Bin81}. Since $B_1 = B$ is a complex space by Step 1, only the induction step has to be performed. 
We consider the canonical morphisms
$$ w: \sO_{B_k} \to \sO_{B_{k-1}}$$
and 
$$ v: \sO_{B_k} \to \sO_B.$$ 
Clearly ${\rm Ker}(w) \cdot {\rm Ker}(v) = 0$. 
Further, ${\rm Ker}(w) $ is a coherent $\sO_B$-module by Step 1.


Thus $B_k$ is a complex space according to  \cite[10.3]{Bin81}.

The holomorphicity of $f_k$ is then clear by construction. 
Indeed we only have to check that that the pull-back $f_k^* \sO_{B_k}$ is mapped into
$\sO_{S_k}$ or equivalently that we have a morphism $\sO_{B_k} \rightarrow (f_k)_* \sO_{S_k}$. Pushing forward the sequence
$$
0 \rightarrow \sI_S^k \rightarrow \sO_X \rightarrow \sO_{S_k} \rightarrow 0
$$
we obtain
$$
0 \rightarrow \varphi_* \sI_S^k \rightarrow \sO_Y \rightarrow (f_k)_* \sO_{S_k}
$$
Thus the morphism $\sO_{B_k} \rightarrow (f_k)_* \sO_{S_k}$ is obtained by composing
$$
\sO_Y/(\sI_B)^k \rightarrow \sO_Y/ \varphi_* \sI_S^k
$$
with $\sO_Y \rightarrow (f_k)_* \sO_{S_k}$.

{\bf Step 5.} 
Let $m \in \N$ be the Cartier index of $\sI=\sO_X(-S)$, so
$$
\sJ := \sO_X(-mS)
$$
is an invertible sheaf. We consider the morphism $f_m: S_m \rightarrow B_m$. The invertible sheaf
$\sJ$ is $f$-ample, so using  \cite[4.6.13(vi)]{EGA}  and the Nakai-Moishezon criterion (which holds for proper, not necessarily reduced schemes), we know that  $\sJ \otimes \sO_{mS}$ is $f_m$-ample.
Thus by Serre vanishing there exists a $k_0 \in \N$ such that
$$
R^i (f_m)_* (\sJ^k \otimes \sO_{mS}) = 0 \qquad \forall \ k \geq k_0, i \geq 1.
$$
The morphism  $\sJ^{k+1} \otimes \sO_X/\sJ \rightarrow J^k  \otimes \sO_X/\sJ$ being zero for all $k \in \N$ we deduce
$$
\sJ^k/\sJ^{k+1} \otimes  \sO_{mS} \simeq \sJ^k \otimes  \sO_{mS} \qquad \forall \ k \in \N.
$$
Thus we can restate the Serre vanishing as
\begin{equation} \label{serre}
R^i (f_m)_* (\sJ^k/\sJ^{k+1} \otimes  \sO_{mS}) = 0 \qquad \forall \ k \geq k_0, i \geq 1.
\end{equation}
Set now
$$
\sK := \sO_X(-mk_0 S),
$$
then we want to show that
$$
R^1 (f_{mk_0})_* (\sK^d/\sK^{d+1} \otimes  \sO_{mk_0S}) = 0 \qquad \forall \ d \geq 1.
$$

Proceeding by induction over $j=1, \ldots, k_0$ we will show more generally that 
$$
R^1 (f_{mj})_* (\sK^d/\sK^{d+1} \otimes  \sO_{mjS}) = 0 \qquad \forall \ d \geq 1.
$$
For the start of the induction $j=1$ recall that
$$
\sK^d/\sK^{d+1}  = \sO_X(-mk_0 d S)/\sO_X(-mk_0 (d+1) S) = \sJ^{k_0 d}/\sJ^{k_0(d+1)}.
$$
We consider the filtration
$$
\sJ^{k_0(d+1)} \subset \ldots \subset \sJ^{k_0 d+1} \subset \sJ^{k_0 d}
$$
with graded pieces isomorphic to $\sJ^e/\sJ^{e+1}$ with $e \geq k_0 d \geq k_0$. Thus we can apply
\eqref{serre} for each graded piece and obtain 
$$
R^1 (f_{m})_* (\sK^d/\sK^{d+1} \otimes  \sO_{mS}) = 0 \qquad \forall \ d \geq 1.
$$
For the induction step $j-1 \rightarrow j$ note first that we have an exact sequence
$$
0 \rightarrow \sO_{mS}(-m(j-1)S) \rightarrow \sO_{mj S} \rightarrow \sO_{m(j-1) S} \rightarrow 0
$$
Since $mj \leq mk_0$ all these sheaves are supported on subschemes of $mk_0 S$, i.e. the scheme
defined by $\sO_X/\sK$.
Since $\sK^d/\sK^{d+1}$ is locally free on  $mk_0 S$, so is its restriction to any subscheme.
Thus the twisted sequence
$$
0 \rightarrow \sO_{mS}(-m(j-1)S) \otimes \sK^d/\sK^{d+1} \rightarrow \sO_{mj S} \otimes \sK^d/\sK^{d+1}\rightarrow \sO_{m(j-1) S} \otimes \sK^d/\sK^{d+1} \rightarrow 0
$$
is still exact and by 
induction we are left to show that
$$
R^1 (f_m)_*  (\sO_{mS}(-m(j-1)S) \otimes \sK^d/\sK^{d+1}) = 0
$$
for all $d \in \N$. Yet
$$
 \sO_{mS}(-m(j-1)S) \otimes \sK^d/\sK^{d+1}
 \simeq \sO_{mS} \otimes \sJ^{dk_0+j-1}/\sJ^{(d+1)k_0}.
$$
Thus we can consider the filtration
$$
\sJ^{(d+1)k_0} \subset \ldots \subset \sJ^{dk_0+j} \subset \sJ^{(d k_0+j-1}
$$
with graded pieces of the form $\sJ^e/\sJ^{e+1}$ with $e \geq d k_0+j-1 \geq k_0$. Thus we can again apply \eqref{serre}.
\end{proof} 

\begin{remark*}
Our proof explains the inaccuracy in \cite{AV84}: while it is possible to replace
the prime ideal $\sO_X(-S)$ by some power $\sO_X(-mS)$, it is not sufficient to verify
the vanishing of the higher direct images for $f: S \rightarrow B$. This has to be done for the extended morphism $f_m: S_m \rightarrow B_m$, a subtlety that \cite[Cor.8.2]{Bin81}
takes into account.
\end{remark*}

\begin{remark} \label{remark-fujiki}
Fujiki \cite[Prop.3]{Fuj75} constructed an example that contradicts the statement \cite[Prop.7.4]{HP16}: more precisely he constructs a complex manifold $M$ containing a projective manifold $F \subset M$
that is a divisor with antiample normal bundle. Thus $F$
can be contracted onto a point by Grauert's theorem, we denote by $\holom{\varphi_M}{M}{M'}$ the contraction. Moreover there exists an affine bundle $\holom{\pi}{X}{M}$ that is not trivial\footnote{The existence of the affine bundle
is guaranted by Fujiki's condition that $H^1(F, N_{F/M}^*) \neq 0$. By Kodaira vanishing this condition never holds
if $F$ is a Fano manifold, as in the setup of Theorem \ref{thmA}.}, but trivial over $F$, i.e. we have
$$
D: = \fibre{\pi}{F} \simeq F \times \C. 
$$
Let $\holom{f}{D}{\C}$ be the projection, then $D \subset X$ satisfies the conditions
\cite[Prop.7.4]{HP16}, but Fujiki shows that it is not possible to contract
$D \subset X$ onto a curve. 

Let us check that this example does not contradict Theorem \ref{thmB}:
let $Z \subset \C$ be any finite set, then we denote 
$$
D^\circ := D \setminus (F \times Z), \qquad X^\circ := X \setminus (F \times Z)
$$
and
$$
f^\circ := f|_{D^\circ}, \qquad \pi^\circ := \pi|_{X^\circ}.
$$
Arguing by contradiction we assume that 
there exists a bimeromorphic morphism
$$
\holom{\varphi^\circ}{X^\circ}{Y^\circ}
$$
onto a normal complex space $Y^\circ$ such that the restriction to $X^\circ \setminus D^\circ$ is an isomorphism and
$\varphi^\circ|_{D^\circ} = f^\circ$.

We follow Fujiki's argument: by a rigidity argument there exists a bimeromorphic
morphism $\holom{\sigma}{Y^\circ}{M'}$ such that $\sigma \circ \varphi^\circ= \varphi_M \circ \pi^\circ$, i.e the bimeromorphic map $\sigma$ contracts $\C \setminus Z$ 
onto the point $m:= \varphi_M(F)$.  Choose now any point $y \in \varphi^\circ(D^\circ) \simeq \C \setminus Z$.
Up to replacing $M'$ by a small analytic neighbourhood $m \in U \subset M'$ (and therefore $M, X^\circ, Y^\circ$ by $\fibre{\varphi_M}{U}, \fibre{\varphi_M \circ \pi}{U}, \fibre{\sigma}{U}$) we can assume without loss of generality that there exists
a Cartier divisor $S_Y \subset Y^\circ$ such that
$$
S_Y \cap \varphi^\circ(D^\circ) = y
$$ 
as a set. The divisor $S_Y$ does not contain the image of the exceptional locus of $\varphi^\circ$, so the strict transform $S \subset X$ coincides with the pull-back of $S_Y$. In particular we have a set-theoretical equality $S \cap D^\circ = F \times y$
and (up to replacing $M'$ by a smaller neighbourhood of $m$), the morphism
$$
\pi^\circ|_S = \pi|_S : S \rightarrow M
$$
is quasi-finite. Note that $S$ is a subsection of $\pi$ over $F$.
Since  $X^\circ \rightarrow M$ is an affine bundle over $M \setminus F$ we can use Fujiki's arithmetic mean construction \cite[p.505]{Fuj75} to construct a section $\bar M \subset X^\circ \subset X$. Yet this contradicts the fact that the affine bundle $X \rightarrow M$ is non-trivial.
\end{remark}

\begin{proof}[Proof of Theorem \ref{thmA}]
Let $f': S \rightarrow B'$ be the Stein factorisation of $f$, then by construction
$f_* \sO_S \simeq \sO_{B'}$ (cf. the proof of \cite[III,Cor.11.5]{Har77}).   
In order to simplify the notation we assume that $f=f'$.

By \cite[7.8]{HP16}, the general fiber of $f$ is $\PP_1$. 
The threefold $X$ being $\Q$-factorial the divisor $S \subset X$ is $\Q$-Cartier, say of index $m \in \N$. Since $f$ contracts an extremal ray it is clear that $\sO_S(-mS)$ is $f$-ample.

Since $X$ is a threefold with terminal singularities, the singular locus ${\rm Sing}(X)$ is finite. Set
$$
Z := f({\rm Sing (X)}) \cup {\rm Sing}(B), 
$$
and
$$
X^{\circ} := X \setminus \fibre{f}{Z}, \ S^{\circ} := S \setminus \fibre{f}{Z}, \ B^\circ := B \setminus Z, \ f^{\circ} = f \vert S^{\circ}.
$$
Then $f^{\circ}: S^\circ \rightarrow B^\circ$ is a $\PP^1$-bundle such that the conormal bundle $N^*_{S^\circ/X^\circ}$
is $f^\circ$-ample.
By the classical theorem of Nakano \cite{Nak71}, \cite{FN71},  there exists a bimeromorphic holomorphic map 
$$
\varphi^{\circ}: X^{\circ} \to Y^{\circ}
$$
to a normal complex space $Y^{\circ}$ such that $\varphi^{\circ} \vert_{S^{\circ}} = f^{\circ}$ and $\varphi^{\circ} $ is the blowup of $Y^\circ$ along the submanifold $B^\circ$. 
Now apply Theorem \ref{thmB}.
\end{proof}


\begin{thebibliography}{Nak71}

\bibitem[AT82]{AT82}
Vincenzo Ancona and Giuseppe Tomassini, \emph{Modifications analytiques},
  Lecture Notes in Mathematics, vol. 943, Springer-Verlag, Berlin, 1982.
  \MR{MR673560 (84g:32022)}

\bibitem[AT84]{AV84}
V.~Ancona and Vo~Van Tan, \emph{On the blowing down problem in {${\bf
  C}$}-analytic geometry}, J. Reine Angew. Math. \textbf{350} (1984), 178--182.
  \MR{743541 (86h:32025)}

\bibitem[Bin81]{Bin81}
J\"{u}rgen Bingener, \emph{On the existence of analytic contractions}, Invent.
  Math. \textbf{64} (1981), no.~1, 25--67. \MR{621769}

\bibitem[DH23]{DH23}
Omprokash Das and Christopher Hacon, \emph{On the minimal model program for
  {K}{\"a}hler 3-folds}, arXiv preprint \textbf{2306.11708} (2023).

\bibitem[FK72]{FK72}
O.~Forster and K.~Knorr, \emph{Relativ-analytische {R}\"{a}ume und die
  {K}oh\"{a}renz von {B}ildgarben}, Invent. Math. \textbf{16} (1972), 113--160.
  \MR{325998}

\bibitem[FN72]{FN71}
Akira Fujiki and Shigeo Nakano, \emph{Supplement to ``{O}n the inverse of
  monoidal transformation''}, Publ. Res. Inst. Math. Sci. \textbf{7} (1971/72),
  637--644. \MR{294712}

\bibitem[Fuj75]{Fuj75}
Akira Fujiki, \emph{On the blowing down of analytic spaces}, Publ. Res. Inst.
  Math. Sci. \textbf{10} (1975), 473--507 (English).

\bibitem[GR79]{GR79}
Hans Grauert and Reinhold Remmert, \emph{Theory of {S}tein spaces}, Grundlehren
  der Mathematischen Wissenschaften, vol. 236, Springer-Verlag, Berlin-New
  York, 1979, Translated from the German by Alan Huckleberry. \MR{580152}

\bibitem[Gro61]{EGA}
A.~Grothendieck, \emph{{\'E}l\'ements de g\'eom\'etrie alg\'ebrique. {II}.
  {\'e}tude globale \'el\'ementaire de quelques classes de morphismes.}, Inst.
  Hautes \'Etudes Sci. Publ. Math. (1961), no.~8, 222. \MR{217084}

\bibitem[Har77]{Har77}
Robin Hartshorne, \emph{Algebraic geometry}, Springer-Verlag, New York, 1977,
  Graduate Texts in Mathematics, No. 52. \MR{MR0463157 (57 \#3116)}

\bibitem[HP16]{HP16}
Andreas H{\"o}ring and Thomas Peternell, \emph{Minimal models for {K}\"ahler
  threefolds}, Invent. Math. \textbf{203} (2016), no.~1, 217--264. \MR{3437871}

\bibitem[Nak71]{Nak71}
Shigeo Nakano, \emph{On the inverse of monoidal transformation}, Publ. Res.
  Inst. Math. Sci. \textbf{6} (1970/71), 483--502. \MR{294710}

\end{thebibliography}

\providecommand{\bysame}{\leavevmode\hbox to3em{\hrulefill}\thinspace}
\providecommand{\MR}{\relax\ifhmode\unskip\space\fi MR }
\providecommand{\MRhref}[2]{%
  \href{http://www.ams.org/mathscinet-getitem?mr=#1}{#2}
}
\providecommand{\href}[2]{#2}

\end{document}